\documentclass[12pt]{article}
\usepackage{cite}
\usepackage{amsmath}
\usepackage{amsthm}
\usepackage{amssymb}
\usepackage[margin=1in]{geometry}
\usepackage[colorlinks=true, linkcolor=blue, citecolor=blue, urlcolor=blue]{hyperref}

\newtheorem{theorem}{Theorem}
\newtheorem{lemma}{Lemma}

\newtheorem{definition}{Definition}
\newtheorem{problem}{Problem}

\begin{document}

\title{A Riemann Boundary Value Problem in a Two-Dimensional Commutative Associative Banach Algebra}
\author{Sergiy Plaksa\thanks{Institute of Mathematics, National Academy of Sciences of Ukraine, Tereshchenkivska~3, 01024 Kyiv, Ukraine. E-mail: \texttt{plaksa62@gmail.com}} \and Roman Pukhtaievych\thanks{Institute of Mathematics, National Academy of Sciences of Ukraine, Tereshchenkivska~3, 01024 Kyiv, Ukraine. E-mail: \texttt{rpukhtaievych@gmail.com}}}
\date{\today}
\maketitle

\begin{abstract}
    We consider a Riemann boundary value problem for monogenic functions in a two-dimensional commutative associative Banach algebra. We prove theorems on the existence of a solution to this problem under different assumptions on the coefficient and free term of the problem, and give an explicit formula for the solution.
\end{abstract}

\medskip
\noindent\textbf{Keywords:} Riemann boundary value problem; monogenic function; commutative Banach algebra; Douglis algebra; dual complex numbers; Cauchy-type integral; Sokhotski-Plemelj formulas.

\medskip
\noindent\textbf{2020 Mathematics Subject Classification:} 30G35 (primary); 30E25, 46J15 (secondary).
\medskip

\section{Introduction}

The Riemann boundary value problem (RBVP) is a classical problem in complex analysis that seeks to find a piecewise analytic function \( \Phi \) in the complement of a Jordan curve \( \gamma \subset \mathbb{C} \) such that its boundary values \( \Phi^+ \) and \( \Phi^- \) from the interior and exterior domains, respectively, satisfy the boundary condition
\[
\Phi^+(x) = G(x)\Phi^-(x) + g(x), \quad x \in \gamma,
\]
where \( G \) is a given coefficient function and \( g \) is a given free term. This problem, first formulated by Riemann in 1851, has applications in mathematical physics, elasticity theory, fluid dynamics, and the theory of singular integral equations.

The classical theory of RBVP was developed over the course of nearly a century through the contributions of many mathematicians. Plemelj (1908, \cite{Ple08}) and Carleman (1922, \cite{Car22}) obtained solutions in particular cases. A major breakthrough came with the work of Gakhov (1937, \cite{Gakhov90}), who obtained the explicit solution of RBVP on smooth curves for coefficient and free term functions from the Hölder classes, under the condition that \( G(x) \neq 0 \) for all \( x \in \gamma \). Gakhov established the fundamental principle that the solvability of RBVP is determined by the index of the problem, defined as
\[
\varkappa := \frac{1}{2\pi i} \int_\gamma d(\ln G(t)).
\]
This index characterizes the winding number of the coefficient \( G \) along the contour and plays a decisive role in determining the dimension of the solution space. Building on this foundation, Gakhov extended his method to handle piecewise continuous boundary value problems by reducing them to the continuous case.

Muskhelishvili (1941, \cite{Muskh41}) developed a modified solution scheme for RBVP on unclosed curves that provided a unified approach applicable to both continuous and piecewise continuous problems, regardless of the contour type.

A natural question that arose from this classical theory was whether the regularity assumptions on the coefficient and free term could be extended beyond the Hölder class. Gerus (1981, \cite{Gerus81}) made 
progress in this direction by extending Gakhov's method to any rectifiable Jordan curve and a Dini-continuous free term $g$ 
but with an additional restriction on the function $G$.
Moreover, Babaev and Salaev (1982, \cite{BaSa82})  employed the analyticity properties of the solution function in the interior domain 
to prove that Gakhov's theorems hold for any rectifiable Jordan curve when both \( G \) and \( g \) are Dini-continuous functions.

In parallel with these developments in classical complex analysis, there has been growing interest in extending boundary value problems to more general algebraic settings. Gilbert and Zeng \cite{GilbertZeng1992}, Blaya, Reyes, and Pe{\~n}a (2005, \cite{AbreuBlaya2005}) considered the RBVP of hyperanalytic functions for the case of data given on rectifiable closed curves in Douglis algebras. 

A natural framework for such generalizations is provided by the theory of monogenic functions in commutative associative Banach algebras. Of our particular interest is the two-dimensional commutative associative Banach algebra \( \mathbb{B} \) over \( \mathbb{C} \) with basis \( \{1, \rho\} \) where \( \rho^2 = 0 \), known as the algebra of dual complex numbers, two-dimensional Douglis algebra, or biharmonic algebra. This algebra has found applications in the study of biharmonic functions and potential theory.

The purpose of this paper is to develop a comprehensive theory of the Riemann boundary value problem for monogenic functions defined on a two-dimensional real linear subspace \( E \) of the Banach algebra \( \mathbb{B} \). We establish fundamental results analogous to the classical theorems of Cauchy, Morera, and Liouville in this setting, and prove Sokhotski-Plemelj formulas for Cauchy-type integrals with densities that are products of monogenic functions and Dini-continuous functions. Our main results provide explicit solutions to jump, homogeneous, and non-homogeneous Riemann boundary value problems in terms of the problem index, and establish solvability conditions when the index is negative. The approach we develop extends the classical methods of Gakhov to the algebra \( \mathbb{B} \), incorporating the techniques of Gerus and Babaev-Salaev.

The paper is organized as follows. In Section 2, we introduce the notation and preliminaries, including the structure of the algebra \( \mathbb{B} \), the definition of monogenic functions, and the class of Dini-continuous functions. Section 3 establishes fundamental properties of monogenic functions, including analogues of classical theorems from complex analysis. In Section 4, we prove the Sokhotski-Plemelj formulas for Cauchy-type integrals in our setting. Finally, Section 5 presents the main results on the Riemann boundary value problem, providing explicit solutions and solvability conditions for various regularity assumptions on the coefficient and free term.

\section{Notation and Preliminaries}

Let \( \mathbb{R} \) and \( \mathbb{C} \) denote the fields of real and complex numbers, respectively.
For any subset \( D \subset \mathbb{C} \), we denote its closure by \( \overline{D} \).
For a point \( z \in \mathbb{C} \), we write \( \Re(z) \) and \( \Im(z) \) for its real and imaginary parts.

Let \( \mathbb{B} \) be the two-dimensional commutative associative Banach algebra over \( \mathbb{C} \) with basis \( \{1, \rho\} \), where \( \rho^2 = 0 \).
Any element \( c \in \mathbb{B} \) has a unique representation \( c = c_1 + c_2 \rho \), where \( c_1, c_2 \in \mathbb{C} \).
The algebra \( \mathbb{B} \) is often referred to as the algebra of dual complex numbers, two-dimensional Douglis algebra, or biharmonic algebra (see, e.g., \cite{KovMel1891, Lorch43, GrPl09}).

An element \( c = c_1 + c_2 \rho \in \mathbb{B} \) is invertible if and only if its complex part \( c_1 \neq 0 \). Its inverse is given by
\[
c^{-1} = \frac{1}{c_1} - \frac{c_2}{c_1^2} \rho.
\]
The set of invertible elements in \( \mathbb{B} \) is denoted by \( \mathbb{B}^{\times} \).

For any \( c \in \mathbb{B}^{\times}\), its logarithm is defined (see \cite[p.~422]{Lorch43}) as
\begin{equation}\label{log-def}
\ln c := \ln c_1 + \frac{c_2}{c_1} \rho,
\end{equation}
where \( \ln c_1 \) is the principal value of the complex logarithm (see, e.g., \cite[p. 229]{GrPl2010}).

We consider a two-dimensional real linear subspace of \( \mathbb{B} \), denoted by \( E \), which is generated by a basis \( \{e_1, e_2\} \).
The one-point compactification of \( E \) is denoted by \( \overline{E} := E \cup \{\infty\} \).
Let
\[
e_1 := a_1 + b_1 \rho, \qquad e_2 := a_2 + b_2 \rho,
\]
where \( a_1, a_2, b_1, b_2 \in \mathbb{C} \) satisfy the condition
\begin{equation} \label{eq:basis}
\det \begin{pmatrix} \Re(a_1) & \Re(a_2) \\ \Im(a_1) & \Im(a_2) \end{pmatrix} \neq 0.
\end{equation}
Thus, any element \( \zeta \in E \) can be uniquely represented as \( \zeta := x e_1 + y e_2 \) for some \( x, y \in \mathbb{R} \),
and as \( \zeta := \xi_1 + \xi_2 \rho \) for some \( \xi_1, \xi_2 \in \mathbb{C} \).  Condition \eqref{eq:basis} ensures that:
\begin{itemize}
    \item[(i)] \( e_1 \) and \( e_2 \) are linearly independent over \( \mathbb{R} \).
    \item[(ii)] Every non-zero element of \( E \) is invertible in \( \mathbb{B} \).
\end{itemize}
The biharmonic basis \( e_1 = 1 \) and \( e_2 = i - \frac{i}{2}\rho \) (see, e.g., \cite{KovMel1891}) is an example that satisfies condition \eqref{eq:basis}.

Let \( \|\cdot\| \) denote the norm in \( \mathbb{B} \) defined by
\[ \|c\| := \sqrt{|c_1|^2 + |c_2|^2} \]
 for \( c = c_1 + c_2 \rho \in \mathbb{B} \),
where \( |\cdot| \) denotes the modulus in \( \mathbb{C} \).
Let \( |\cdot|\) also denote the modulus on \( E \) induced by the Euclidean norm in \( \mathbb{R}^2 \), i.e., for \( \zeta = x e_1 + y e_2 \in E \),
\[|\zeta| := \sqrt{x^2 + y^2}.\]
Using the Cauchy-Schwarz inequality, one can verify that
\[ \|\zeta\| \leq c |\zeta| \quad \text{for all } \zeta \in E, \]
where \( c := \sqrt{\|e_1\|^2 + \|e_2\|^2} \) is a constant depending on the choice of the basis \( \{e_1, e_2\} \).

We note that \( E \) is topologically isomorphic to \( \mathbb{C} \) as a real linear space.
The isomorphism between \( E \) and \( \mathbb{C} \) allows us to extend topological concepts from the complex plane to \( E \).
For a Jordan curve \( \gamma \subset E \), \( \overline{E} \setminus \gamma \) consists of a bounded interior domain \( D^+_\gamma \) and an unbounded exterior domain \( D^-_\gamma \).

We now define monogenic functions in this setting.
\begin{definition}
Let \( \Omega \) be a domain in \( \overline{E} \). A function \( \Phi: \Omega \to \mathbb{B} \) is called \textbf{monogenic} in \( \Omega \)
if it is continuous and at every point \( \zeta \in \Omega \) the limit
\[
\Phi'(\zeta) := \lim_{\substack{h \to 0,\, h \in E}} (\Phi(\zeta + h) - \Phi(\zeta))h^{-1}
\]
exists.
\end{definition}

Let us introduce the norm of the function \( \Phi: \Omega \to \mathbb{B} \) as
\[
\|\Phi\|_{\Omega} := \sup_{\zeta \in \Omega} \|\Phi(\zeta)\|,
\]
where \( \Omega \subset \overline{E} \) is the domain of definition of \( \Phi \).

Let \( \gamma \) be a rectifiable curve in \( E \). For each \( \tau \in \gamma \),
\[ \theta_\tau(\epsilon) := \mathrm{mes}\{ t \in \gamma : |t - \tau| \leq \epsilon \}, \] where \( \mathrm{mes} \) is the linear Lebesgue measure on \( \gamma \).
Let \( \mathcal{D}_\gamma \) be the class of functions \( g: \gamma \to \mathbb{B} \) satisfying the Dini condition:
\begin{equation}\label{dini}
    \sup_{\tau\in\gamma}\int_0^1 \frac{\omega_g(\eta)}{\eta} \, d \theta_\tau(\eta) < \infty,
\end{equation}
where the integral is understood in the sense of the upper Darboux integral and
\[ \omega_g(\epsilon) := \sup_{|t_1 - t_2| \leq \epsilon} \|g(t_1) - g(t_2)\| \]
 is the modulus of continuity of \( g \).

For a Jordan curve \( \gamma \in E \), we define the following classes of functions related to the domains \( D^\pm_\gamma \):
\begin{align*}
    \mathcal{H}^+_\gamma &:= \{f: \overline{D^+_\gamma} \to \mathbb{B} \mid f \text{ is monogenic in } D^+_\gamma \text{ and continuous on } \overline{D^+_\gamma}\}, \\
    \mathcal{H}^-_\gamma &:= \{f: \overline{D^-_\gamma} \to \mathbb{B} \mid f \text{ is monogenic in } D^-_\gamma \text{, continuous on } \overline{D^-_\gamma}, \\
    &\qquad \text{ and bounded at } \infty\}, \\
    \mathcal{H}_\gamma &:= \{f:\gamma\to\mathbb{B} \mid f = f^+|_\gamma + f^-|_\gamma \text{ for some } f^\pm\in\mathcal{H}^\pm_\gamma\}.
\end{align*}

For a continuous function \( \psi: \gamma \to \mathbb{B} \), the Cauchy-type integral is defined as
\begin{equation}\label{ITC}
    \widetilde{\psi}(\zeta) := \frac{1}{2\pi i} \int_{\gamma} \psi(\tau)(\tau - \zeta)^{-1}\, d\tau, \quad \zeta \in \overline{E} \setminus \gamma.
\end{equation}
The function \( \widetilde{\psi} \) is monogenic in \( \overline{E} \setminus \gamma \). We denote its limiting values on \( \gamma \) from \( D^\pm_\gamma \) by
\[
\widetilde{\psi}^+(\zeta_0) := \lim_{\zeta \to \zeta_0, \zeta\in D^+_{\gamma}}\widetilde{\psi}(\zeta), \quad
\widetilde{\psi}^-(\zeta_0) := \lim_{\zeta \to \zeta_0, \zeta\in D^-_{\gamma}}\widetilde{\psi}(\zeta), \quad \zeta_0 \in \gamma.
\]

\section{Preliminaries on Monogenic Functions}

In this section, we collect several fundamental theorems for monogenic functions defined on \( \overline{E} \).
These results are analogues of classical theorems from complex analysis
and are essential for the subsequent analysis of the Riemann boundary value problem.

The curvilinear analogues of the Cauchy integral theorem, the Morera theorem, and the Cauchy integral formula have been established
in \cite{Shp15} for the case when \( e_1 = 1 \).
For the specific case of the algebra \( \mathbb{B} \) when \( E \) is defined on a biharmonic basis mentioned earlier in this paper,
these results, along with Taylor and Laurent series expansions, were developed in \cite{GrPl09}.

We state these key theorems below for completeness and note that
they do not follow neither from \cite{Shp15} nor from \cite{GrPl09} directly,
since they are formulated for monogenic functions defined on domains in \( E \) for arbitrary bases \( \{e_1, e_2\} \) satisfying condition \eqref{eq:basis}.
Despite the differences in the settings, the proofs of these theorems can be adapted from the original sources.

\begin{theorem}[Cauchy's Integral Theorem]\label{T:CauchyIntegralTheorem}
    Let a domain \( \Omega \subset E \) and \( \partial \Omega \) be a rectifiable Jordan curve.
    If a function \( \Phi: \overline{\Omega} \to \mathbb{B} \) is monogenic in a domain \( \Omega \) and continuous in \( \overline{\Omega} \),
    then
    \[
    \int_{\partial \Omega} \Phi(\zeta) \, d\zeta = 0.
    \]
\end{theorem}

\begin{theorem}[Morera's Theorem]\label{T:Morera}
    Let \( \Omega \) be a domain in \( E \), and let a function \( \Phi: \Omega \to \mathbb{B} \) be continuous in \( \Omega \).
    If for every triangle \( T \) such that \( \overline{T} \subset \Omega \), the equality
    \[
    \int_{\partial T} \Phi(\zeta) \, d\zeta = 0
    \]
    holds, then \( \Phi \) is monogenic in \( \Omega \).
\end{theorem}

We note that a triangle \( T \) formed by three vertices \( \zeta_1, \zeta_2, \zeta_3 \in E \) is a domain in \( E \) defined as
\[
T := \left\{ \zeta = \sum_{k=1}^3 \lambda_k \zeta_k \mid \lambda_k > 0 \text{ for } k=1,2,3, \text{ and } \sum_{k=1}^3 \lambda_k = 1 \right\}.
\]
The boundary of the triangle \( T \) is denoted by \( \partial T \) and can be expressed as follows:
\[
\partial T := \left\{ \zeta = \sum_{k=1}^3 \lambda_k \zeta_k \mid \lambda_k \geq 0 \text{ for } k=1,2,3, \sum_{k=1}^3 \lambda_k = 1 , \text{ and }  \prod_{k=1}^3 \lambda_k = 0 \right\}.
\]

\begin{theorem}[Cauchy's Integral Formula]\label{T:CauchyIntegralFormula}
Let a domain \( \Omega \subset E \) and \( \partial \Omega \) be a rectifiable Jordan curve.
    If a function \( \Phi: \overline{\Omega} \to \mathbb{B} \) is monogenic in a domain \( \Omega \) and continuous in \( \overline{\Omega} \),
    then for any point \( \zeta_0 \in \Omega \), the following equality holds:
    \[
    \Phi(\zeta_0) = \frac{1}{2\pi i} \int_{\partial\Omega} {\Phi(\zeta)}(\zeta - \zeta_0)^{-1} \, d\zeta.
    \]
\end{theorem}

A direct consequence of the Cauchy integral formula is the existence of a convergent power series representation for monogenic functions,
analogous to the Taylor series in complex analysis (see, e.g.~\cite[p.~171]{Plaksa_2023}).

\begin{theorem}[Taylor's Theorem]\label{T:Taylor}
    Let \( \Phi: \Omega \to \mathbb{B} \) be a monogenic function in a domain \( \Omega \subset \overline{E} \). For any \( \zeta_0 \in \Omega \), \( \Phi \) can be expanded into a power series
    \[
    \Phi(\zeta) = \sum_{n=0}^{\infty} c_n (\zeta - \zeta_0)^n, \quad \text{where} \quad c_n
        = \frac{\Phi^{(n)}(\zeta_0)}{n!} = \frac{1}{2\pi i} \int_{\gamma} {\Phi(\tau)}(\tau - \zeta_0)^{-n-1} \, d\tau,
    \]
    where \( \gamma \subset \Omega \) is any rectifiable Jordan curve such that \( \zeta_0 \in D^+_\gamma \) and \( \overline{D^+_\gamma} \subset \Omega \).
\end{theorem}

Using the Taylor series representation, one can establish the following versions of Liouville's theorem for monogenic functions on \( E \), similar to the classical results (cf. \cite[p.~108]{ShabatP1}).

\begin{theorem}[Liouville's Theorem]\label{T:Liouville}
    If a function \( \Phi: \overline{E} \to \mathbb{B} \) is monogenic and bounded at every point of \( \overline{E} \), then it must be a constant.
\end{theorem}

\begin{theorem}[Extended Liouville's Theorem]\label{T:Liouville-Extended}
    Let \( \Phi: \overline{E} \to \mathbb{B} \) be a monogenic function on \( \overline{E} \). If \( \Phi \) has a pole of order at most \( m \) at infinity, then \( \Phi \) is a polynomial of degree at most \( m \).
\end{theorem}

Finally, we state a principle of monogenic continuation (analogous to the principle of analytic continuation in complex analysis), which is crucial for solving boundary value problems by piecing together solutions defined on adjacent domains (cf. \cite[Sec. 19]{Gakhov90}).

\begin{lemma}[Principle of Monogenic Continuation]\label{L:MonogenicCont}
    Let \( \Omega_1 \) and \( \Omega_2 \) be two domains in \( \overline{E} \) whose boundaries share a common open rectifiable Jordan arc \( \gamma \). Let \( \Phi_1: \Omega_1 \cup \gamma \to \mathbb{B} \) and \( \Phi_2: \Omega_2 \cup \gamma \to \mathbb{B} \) be continuous functions that are monogenic in \( \Omega_1 \) and \( \Omega_2 \), respectively. If \( \Phi_1(\zeta) = \Phi_2(\zeta) \) for all \( \zeta \in \gamma \), then the function
    \[
    \Phi(\zeta) :=
    \begin{cases}
    \Phi_1(\zeta), & \zeta \in \Omega_1 \cup \gamma, \\
    \Phi_2(\zeta), & \zeta \in \Omega_2,
    \end{cases}
    \]
    is monogenic in the domain \( \Omega := \Omega_1 \cup \gamma \cup \Omega_2 \).
\end{lemma}

Next, we will also need the result which can be proved in a similar way to Theorem 2.4 from \cite{MePl08}: every monogenic in a domain \( \Omega \subset E \) function \( \Phi: \Omega \to \mathbb{B} \) can be expressed in the form
\[ \Phi(\zeta) = \frac{1}{2\pi i} \int_{\Gamma} (t - \zeta) ^{-1} F(t) \, dt + \Phi_0(\zeta)\rho ,\]
where \( F \) is a holomorphic function in a domain \( \Omega_{\mathbb{C}}:=\{\xi_1 \in \mathbb{C} \mid \forall \zeta=\xi_1+\xi_2\rho\in \Omega \}\),
\( \Gamma \) is a rectifiable Jordan curve in \( \Omega_{\mathbb{C}} \) such that \( \xi_1 \in D^+_\Gamma \) and \( \overline{D^+_\Gamma} \subset \Omega_{\mathbb{C}} \),
and \( \Phi_0(\zeta)\rho \) is a monogenic function for \( \zeta \in \Omega \).

The proof of the next theorem follows the classical approach used in complex analysis, adapted to the setting of monogenic functions on \( E \): we consider a small neighborhood around each zero \( \zeta_k \) and apply the argument principle to count the zeros.
The details involve constructing appropriate contours and using the properties of monogenic functions, similar to the complex case.

Let \( N_{F} \) denote the number of zeros (counted with multiplicities) of the holomorphic function \( F \) in the domain \( \Omega_{\mathbb{C}} \).

\begin{theorem}[Logarithmic Residue Theorem]\label{T:LogResidue}
    Let \( \gamma \subset E  \) be a rectifiable Jordan curve. Let \( D_\gamma \subset \overline{E} \) be \( D^+_\gamma \) or \( D^-_\gamma \).
    Let  \( D\)  be a domain such that \( \overline{D_\gamma} \subset D \).
    Let \( \Phi: D \to \mathbb{B} \) be a function monogenic in \( D\), continuous on \( \overline{D_\gamma} \),
    and \( \Phi(\zeta) \in \mathbb{B}^{\times} \) for all \( \zeta \in \gamma \). Then the following equality holds:
    \[
    \frac{1}{2\pi i} \int_{\gamma} {\Phi'(\zeta)}(\Phi(\zeta))^{-1} \, d\zeta =
    \begin{cases}
    N_{F}, & D_\gamma \text{ is } D^+_\gamma, \\
    - N_{F}, & D_\gamma \text{ is } D^-_\gamma.
    \end{cases}
    \]
\end{theorem}

This is a particular case of the more general argument principle for monogenic functions which takes into account the zeros of the holomorphic function \( F \) in the representation of \( \Phi \)
and singularities of \( \Phi \) itself.

\section{Sokhotski-Plemelj Formulas}
The following theorem is an analogue of the Sokhotski-Plemelj formulas for Cauchy-type integrals with densities that are products of monogenic functions and Dini-continuous functions.

\begin{theorem}\label{T:SokhPlem_H+}
    Let \( \gamma \subset E \) be a rectifiable Jordan curve, let \( g \in \mathcal{D}_\gamma \), and let a function \( H^+ \in \mathcal{H}^+_\gamma \). Then for the function
    \[
    \psi(\tau) := H^+(\tau)g(\tau) \quad \forall \tau \in \gamma,
    \]
    the Cauchy-type integral \( \widetilde{\psi} \) has limiting values \( \widetilde{\psi}^\pm(\tau) \) at every point \( \tau \in  \gamma \) from the domains $D^\pm_\gamma$, and the following equality holds:
    \[
    \psi(\tau) = \widetilde{\psi}^+(\tau) - \widetilde{\psi}^-(\tau) \quad \forall \tau \in \gamma.
    \]
\end{theorem}

The following analogous result holds for the case when \( H^- \in \mathcal{H}^-_\gamma \).

\begin{theorem}\label{T:SokhPlem_H-}
    Let \( \gamma \subset E \) be a rectifiable Jordan curve, let \( g \in \mathcal{D}_\gamma \),
    and let \( H^- \in \mathcal{H}^-_\gamma \). Then for the function
    \[
    \psi(\tau) := H^-(\tau)g(\tau) \quad \forall \tau \in \gamma,
    \]
    the Cauchy-type integral \( \widetilde{\psi} \) has limiting values \( \widetilde{\psi}^\pm(\tau) \) at every point \( \tau \in \gamma \) from the domains \( D^\pm_\gamma \), and the following equality holds:
    \[
    \psi(\tau) = \widetilde{\psi}^+(\tau) - \widetilde{\psi}^-(\tau) \quad \forall \tau \in \gamma.
    \]
\end{theorem}

\section{Riemann Boundary Value Problems for Monogenic Functions}
In this section, we consider the Riemann boundary value problem for monogenic functions defined on the real linear subspace \( E \) of the algebra \( \mathbb{B} \).
\begin{problem}[Riemann Boundary Value Problem]\label{Pb:BVP-NHom}
    Let \( \gamma \subset E \) be a rectifiable Jordan curve and let \( G: \gamma \to \mathbb{B} \) and \( g: \gamma \to \mathbb{B} \) be given functions.
    Find a pair of functions \( (\Phi^+, \Phi^-) \in \mathcal{H}^+_\gamma \times \mathcal{H}^-_\gamma \) satisfying the boundary condition
    \begin{equation} \label{bc-nhom}
        \Phi^+(\tau) = G(\tau)\Phi^-(\tau) + g(\tau) \quad \forall \tau \in \gamma.
    \end{equation}
\end{problem}
Here, the function \( G \) is called the \textbf{coefficient} of the problem, and the function \( g \) is called
the \textbf{free term} or \textbf{non-homogeneous term} of the problem. If \( g \equiv 0 \),
then Problem \ref{Pb:BVP-NHom} is called a \textbf{homogeneous Riemann boundary value problem}.

For a rectifiable Jordan curve \( \gamma \subset E \) and function \( G: \gamma \to \mathbb{B}^{\times} \) from Problem \ref{Pb:BVP-NHom},
we define the \textbf{index} of the problem as
\[
\varkappa_{G, \gamma} := \frac{1}{2\pi i}\int_{\gamma} d(\ln G(t)).
\]

Below we state four theorems that provide the solution to Problem \ref{Pb:BVP-NHom} under different settings for \( g \) and \( G \) in terms of the index  \( \varkappa := \varkappa_{G, \gamma} \).

\begin{theorem}\label{Th:Jump}
    Let \( \gamma \subset E \) be a rectifiable Jordan curve,
    a function \( G\equiv1 \),
    and a function \( g \in \mathcal{H}_{\gamma} \) in Problem~\ref{Pb:BVP-NHom}.
    Then the solution of Problem~\ref{Pb:BVP-NHom} is given by
    \[
    \Phi^\pm(\zeta) = \widetilde{g}(\zeta) + c \quad \forall \zeta \in D^\pm_{\gamma},
    \]
    where \( \widetilde{g} \) is the Cauchy-type integral of \( g \) defined by \eqref{ITC}, and \( c \) is an arbitrary constant in \( \mathbb{B} \).
\end{theorem}

In the following theorem, we consider the homogeneous Riemann boundary value problem, i.e., Problem~\ref{Pb:BVP-NHom} with \( g \equiv 0 \), and provide the solution in terms of the index \( \varkappa \).

\begin{theorem}\label{Th:HomBVP}
    Let \( \gamma \subset E \) be a rectifiable Jordan curve and \( 0 \in D^+_\gamma \). Let a function \( \ln(\tau^{-\varkappa} G(\tau)) \in \mathcal{H}_\gamma \) with \( G(\tau) \in \mathbb{B}^{\times} \) for all \( \tau \in \gamma \), and \( g \equiv 0 \) in Problem~\ref{Pb:BVP-NHom}.

    Define the function \( X_0: \overline{E} \setminus \gamma \to \mathbb{B} \) by
    \[
    X_0(\zeta) := \exp\left( \frac{1}{2\pi i} \int_\gamma (\tau-\zeta)^{-1} \ln (\tau^{-\varkappa}G(\tau))\, d\tau \right),
    \]
    and set
    \[
        X(\zeta) :=
        \begin{cases}
            X_0(\zeta), & \zeta \in D^+_\gamma, \\
            \zeta^{-\varkappa}X_0(\zeta), & \zeta \in D^-_\gamma.
        \end{cases}
    \]
    Let \( X^\pm \) denote the limiting values of \( X \) on \( \gamma \) from \( D^\pm_\gamma \).

    Then Problem~\ref{Pb:BVP-NHom} has the following solutions:
    \begin{itemize}
        \item[(i)] If \( \varkappa \geq 0 \), the general solution is
        \[
        \Phi^\pm(\zeta) = X^\pm(\zeta) P_{\varkappa}(\zeta) \quad \forall \zeta \in D^\pm_{\gamma},
        \]
        where \( P_{\varkappa} \) is an arbitrary polynomial of degree at most \( \varkappa \).
        \item[(ii)] If \( \varkappa < 0 \), the problem has only the trivial solution \( \Phi^\pm \equiv 0 \).
    \end{itemize}
\end{theorem}

In the following two theorems, we consider the non-homogeneous Riemann boundary value problem under different settings of \( G\) and \( g \), and provide their solutions in terms of the index \( \varkappa \) and the solvability conditions.

\begin{theorem}\label{Th:NHBVP-1}
    Let \( \gamma \subset E \) be a rectifiable Jordan curve and \( 0 \in D^+_\gamma \). Let a function \( G \) belong to \( \mathcal{D}_\gamma \) with \( G(\tau) \in \mathbb{B}^{\times} \) for all \( \tau \in \gamma \), and a function \( g \in \mathcal{H}_{\gamma} \). Let \( X_0 \) and \( X \) by defined as in Theorem~\ref{Th:HomBVP}. Let
    \[
    \psi(\tau) := g(\tau)(X^+(\tau))^{-1} \quad \forall \tau \in \gamma.
    \]

    Then Problem~\ref{Pb:BVP-NHom} is solvable if and only if one of the following conditions holds:
    \begin{itemize}
        \item[(i)] \( \varkappa \geq 0 \);
        \item[(ii)] \( \varkappa < 0 \) and the solvability conditions
                \begin{equation}\label{eq:solvability}
                \int_{\gamma} \psi(\tau) \tau^{s-1} \, d\tau = 0, \quad s = 1, 2, \dots, -\varkappa,
                \end{equation}
                are satisfied.
    \end{itemize}
    Moreover, the general solution of Problem~\ref{Pb:BVP-NHom} is given by
    \begin{equation}\label{eq:general_solution}
    \Phi^\pm(\zeta) = X^\pm(\zeta) \left( \widetilde{\psi}(\zeta) + P_{\varkappa}(\zeta) \right) \quad \forall \zeta \in D^\pm_{\gamma},
    \end{equation}
    where \( \widetilde{\psi} \) is the Cauchy-type integral \eqref{ITC} of \( \psi \), and \( P_{\varkappa} \) is an arbitrary polynomial of degree at most \( \varkappa \) when \( \varkappa \geq 0 \), and \( P_{\varkappa} \equiv 0 \) when \( \varkappa < 0 \).
\end{theorem}

\begin{theorem}\label{Th:NHBVP-2}
    Let \( \gamma \subset E \) be a rectifiable Jordan curve and \( 0 \in D^+_\gamma \). Let a function \( \ln(\tau^{-\varkappa} G(\tau)) \in \mathcal{H}_\gamma \) with \( G(\tau) \in \mathbb{B}^{\times} \) for all \( \tau \in \gamma \), and a function \( g \in \mathcal{D}_{\gamma} \). Let \( X_0 \), \( X \), and \(\psi(\tau)\) be defined as in Theorem~\ref{Th:HomBVP}.

    Then Problem~\ref{Pb:BVP-NHom} is solvable if and only if \( \varkappa \geq 0 \) or \( \varkappa < 0 \) and the solvability conditions \eqref{eq:solvability} are satisfied. Moreover, the general solution of Problem~\ref{Pb:BVP-NHom} is given by \eqref{eq:general_solution}.
\end{theorem}

\section*{Acknowledgments}

This work was supported by a grant from the Simons Foundation (SFI-PD-Ukraine-00014586, R.P.P.) for the second author. The first author was supported by the National Research Foundation of Ukraine, Project number 2025.07/0014, Project name: ``Modern problems of Mathematical Analysis and Geometric Function Theory''.

\bibliographystyle{alpha} 

\begin{thebibliography}{BRP05}

\bibitem[BRP05]{AbreuBlaya2005}
R.~A. Blaya, B.~Juan Reyes, and D.P. Pe{\~n}a.
\newblock Riemann boundary value problem for hyperanalytic functions.
\newblock {\em Int. J. Math. Math. Sci.}, 2005(17):2821--2840, 2005.

\bibitem[BS82]{BaSa82}
A.~A. Babaev and V.~V. Salaev.
\newblock Boundary value problems and singular equations on a rectifiable contour.
\newblock {\em Mat. Zametki}, 31(4):571--580, 654, 1982.

\bibitem[Car22]{Car22}
T.~Carleman.
\newblock Sur la r{\'e}solution de c{\'e}rtaines {\'e}quations int{\'e}grales.
\newblock {\em Ark. Mat. Astron. Fys.}, 16(26), 1922.

\bibitem[Gak90]{Gakhov90}
F.~D. Gakhov.
\newblock {\em Boundary value problems}.
\newblock Dover Publications, Inc., New York, 1990.

\bibitem[Ger81]{Gerus81}
O.~F. Gerus.
\newblock A method of integral equations and the riemann boundary problem.
\newblock {\em Ukr. Math. J.}, 33(3):293--296, May 1981.

\bibitem[GP09]{GrPl09}
S.~V. Grishchuk and S.~A. Plaksa.
\newblock Monogenic functions in the biharmonic plane.
\newblock {\em Reports of the National Academy of Sciences of Ukraine}, 7(12):13--20, 2009.

\bibitem[GP10]{GrPl2010}
S.~Grishchuk and S.~Plaksa.
\newblock On logarithmic residue of monogenic functions of biharmonic variable.
\newblock In {\em Zb. Pr. Inst. Mat. NAN Ukr.}, volume~7, pages 227--234. Institute of Mathematics of NAS of Ukraine, Kyiv, Ukraine, 2010.

\bibitem[GY92]{GilbertZeng1992}
R.~P. Gilbert and Zeng Yuesheng.
\newblock Hyperanalytic riemann boundary value problems on rectifiable closed curves.
\newblock {\em Complex Var. Theory Appl.}, 20(1--4):277--288, 1992.

\bibitem[KM81]{KovMel1891}
V.~F. Kovalyov and I.~P. Melnychenko.
\newblock Biharmonic functions on the biharmonic plane.
\newblock {\em Dokl. Akad. Nauk Ukr. SSR. Series A}, 8:25--27, 1981.

\bibitem[Lor43]{Lorch43}
E.~R. Lorch.
\newblock The theory of analytic functions in normed {A}belian vector rings.
\newblock {\em Trans. Amer. Math. Soc.}, 54:414--425, 1943.

\bibitem[MP08]{MePl08}
I.~P. Mel'nichenko and S.~A. Plaksa.
\newblock {\em Commutative algebras and spatial potential fields}.
\newblock Kyiv: Inst. Math. NAS Ukraine, 2008.

\bibitem[Mus41]{Muskh41}
N.~I. Muskhelishvili.
\newblock Application of integrals of cauchy type to a class of singular integral equations.
\newblock {\em Trans. Inst. Math. Tbilissi}, 10(1), 1941.

\bibitem[Ple08]{Ple08}
J.~Plemelj.
\newblock Ein {E}rg\"anzungssatz zur {C}auchyschen {I}ntegraldarstellung analytischer {F}unktionen, {R}andwerte betreffend.
\newblock {\em Monatsh. Math. Phys.}, 19(1):205--210, 1908.

\bibitem[PS23]{Plaksa_2023}
S.~A. Plaksa and V.~S. Shpakivskyi.
\newblock {\em Monogenic Functions in Spaces with Commutative Multiplication and Applications}.
\newblock Springer Nature Switzerland, 2023.

\bibitem[Sha85]{ShabatP1}
B.~V. Shabat.
\newblock {\em Introduction to complex analysis. Part 1: Functions of one variable}.
\newblock AMS, Providence, 1985.

\bibitem[Shp15]{Shp15}
V.S. Shpakivskyi.
\newblock Integral theorems for monogenic functions in commutative algebras.
\newblock In {\em Zb. Pr. Inst. Mat. NAN Ukr.}, volume 12(4), pages 313--328. Institute of Mathematics of NAS of Ukraine, Kyiv, Ukraine, 2015.

\end{thebibliography}

\end{document}